\documentclass[12pt,reqno,intlimits]{article}

\usepackage{amssymb,amsthm,amsmath,epsfig,graphicx}

\theoremstyle{remark}

\newcommand{\al}{\alpha}
\newcommand{\sbs}{\subset}
\newcommand{\bt}{\beta}
\newcommand{\lb}{\lambda}

\newcommand{\intl}{\int\limits}

\newcommand{\bR}{\Bbb{R}}
\newcommand{\bZ}{\Bbb{Z}}
\newcommand{\Om}{\Omega}
\newcommand{\wt}{\widetilde}

\newcommand{\ve}{\varepsilon}

\newcommand{\const}{\operatorname{const}}

\newcommand{\ol}{\overline}

\begin{document}
\begin{center}
{\bf TWO-WEIGHTED NORM INEQUALITIES FOR THE DOUBLE HARDY
TRANSFORMS AND STRONG FRACTIONAL MAXIMAL FUNCTIONS IN VARIABLE
EXPONENT LEBESGUE SPACES}

\end{center}

\vskip+0.5cm

\begin{center}
{\bf Vakhtang Kokilashvili and Alexander Meskhi}
\end{center}
\vskip+0.5cm

{\bf Abstract.} Two--weight norm estimates for the
double Hardy transforms and strong fractional maximal functions
are established in variable exponent Lebesgue spaces. Derived conditions are simultaneously necessary and sufficient in the case when the exponent of the right--hand  side space is constant.

\vskip+0.5cm


\vskip+0.5cm

 \vskip+1cm

\vskip+0.5cm

\section{INTRODUCTION}

Our  goal is to establish two--weight criteria for the double
Hardy transforms and  strong fractional maximal functions  in
variable exponent Lebesgue spaces. Our interest in these  problems
is stipulated by the following circumstances: by the need in
various applications to the boundary value problems in PDE and the
fact that  strong maximal operator, unlike of the Hardy-Littlewood
maximal function is bounded in $L^{p(\cdot)}$ space if and only if
$p(x)\equiv\text{const}$ (see \cite{Kop1}). As we shall see below
the similar phenomenon occurs also for strong fractional maximal
operators.

\vskip+0.1cm
Let us recall some well--known  results for the classical Lebesgue spaces (see e.g, \cite{KuMaPe}, \cite{KuPe}).

\vskip+0.2cm

The celebrated classical Hardy inequality states:

\vskip+0.2cm

{\bf Theorem A.} {\em Let $p$ be constant satisfying the condition
$1<p<\infty$ and let $f$ be a measurable, nonnegative function in
$(0,\infty)$. Then}
$$ \bigg(\int_0^\infty\bigg(\frac{1}{x}\int_0^x\ f(y) dy
\bigg)^{\frac{1}{p}}\bigg)\leq \frac{p}{p-1}\bigg(\int_0^\infty
f^p(x)dx\bigg)^{\frac{1}{p}}. $$

\vskip+0.2cm

Two-weighted boundedness criteria for the Hardy transform
$$ (\mathcal{H}_{1}f)(x)=\int\limits_{0}^{x}f(y)dy, $$
reads as follows:

\vskip+0.2cm

{\bf Theorem B.} {\em Let $p$ and $q$ be constants satisfying the
condition $1<p\leq q<\infty $. Suppose that $u$ and $v$ are weight
functions  on $\mathbb{R}_{+}.$ Then each of the following
conditions are necessary and sufficient for the inequality
\begin{equation}
 \left( \int\limits_{0}^{\infty }\left( \int\limits_{0}^{x}f(t)dt\right)
 ^{q}v(x)dx\right) ^{\frac{1}{q}}\leq C\left( \int\limits_{0}^{\infty
 }f^{p}(x)w(x)dx\right) ^{\frac{1}{p}}
 \end{equation}
 to hold for all positive and measurable functions on $\mathbb{R}_{+}:$

 a) The Muckenhoupt condition,

$$  A_{M}:=\sup_{x>0} \left( \int\limits_{x}^{\infty }v(t)dt\right)^{
\frac{1}{q}}\left( \int\limits_{0}^{x}w(t)^{1-p^{\prime }}dt\right) ^{\frac{1}{p^{\prime }}}<\infty .  $$

Moreover, the best constant $C$ in $(1.1)$ can be estimated as
follows:
$$
 A_{M}\leq C\leq \left( 1+\frac{q}{p^{\prime }}\right) ^{\frac{1}{q}}\left( 1+ \frac{p^{\prime }}{q}\right) ^{\frac{1}{p^{\prime }}}A_{M}.
$$

 b) The condition of L. E Persson and V. D. Stepanov,
$$
 A_{PS}:= \sup_{x>0} W(x)^{-\frac{1}{p}}\left(
 \int\limits_{0}^{x}v(t)W(t)^{q}dt\right) ^{\frac{1}{q}}<\infty ,\;\;\;
 W(x):=\int\limits_{0}^{x}w(t)^{1-p^{\prime }}dt.
$$
Moreover, the best constant $C$ in $(1.1)$ satisfies the following
estimates:
$$
 A_{PS}\leq C\leq p^{\prime }A_{PS}.$$}
\vskip+0.2cm

In 1984 E.  Sawyer \cite{Saw} found a characterization of
two-weight inequality in terms of three independent conditions for
the double Hardy transform
$$
(\mathcal{H}_{2}f)(x,y)=\int\limits_{0}^{x}\int\limits_{0}^{y}f(t,\tau)dtd\tau.
$$

The following statements gives two-weight criteria in terms of
just one condition when the weight on the right-hand side is a
product of two weights of single variables (see \cite{Mes},
\cite{Ush}, \cite{KokMesPer}, Ch.1):

\vskip+0.2cm

{\bf Theorem C.} {\em Let $p$ and $q$ be constants such that
$1<p\leq q<\infty$ and let $w(x,y)=w_{1}(x)w_{2}(y)$. Then the
operator $\mathcal{H}_{2}$ is bounded from $L^{p}_{w}$ to
$L^{q}_{v}$ $(1<p\leq q<\infty)$ if and only if the Muckenhoupt's
type condition
$$ \sup_{y_{1},y_{2}>0} \left( \int\limits_{y_{1}}^{\infty
 }\int\limits_{y_{2}}^{\infty }v(x_{1},x_{2})dx_{1}dx_{2}\right) ^{\frac{1}{q}%
 }\left(
 \int\limits_{0}^{y_{1}}\int\limits_{0}^{y_{2}}w(x_{1},x_{2})^{1-p^{\prime
 }}dx_{1}dx_{2}\right) ^{\frac{1}{p^{\prime }}}:=A_{1}<\infty
 $$
 is fulfilled.}

\vskip+0.2cm

It should be emphasized that from the results regarding
the two--weight problem derived in this paper, as a corollary, we
deduce trace inequality criteria for the double Hardy transform
when the exponent of the initial Lebesgue space is a constant.
Another remarkable corollary is that there exists a variable
exponent $p(x)$ for which the double average operator is bounded
in $L^{p(\cdot)}$.

In the paper \cite{KoMeJKMS} the authors  established trace
inequality criteria for the strong fractional maximal operator

$$
(M_{\alpha,\beta}f)(x,y):= \sup_{I\times J \ni (x,y)}
\frac{1}{|I|^{1-\alpha}|J|^{1-\beta}}
  \iint\limits_{I\times J} |f(t,\tau)| dt d\tau, \;\;\; 0<\alpha,
  \beta<1,
$$
in constant exponent Lebesgue spaces. In particular, the next
statement holds:

\vskip+0.2cm

 {\bf Theorem D (\cite{KoMeJKMS}).} {\em Let $p$, $q$, $\alpha$ and $\beta$ be constants satisfying the conditions $1<p<q<\infty$ and let
$0<\alpha,\beta<1/p$. Then the following statements are
equivalent:

\rm{(i)} $M_{\alpha,\beta}$ {\em is bounded from}
$L^p({\Bbb{R}}^{2})$ {\em to} $L^q_v({\Bbb{R}}^2)$;

\rm{(ii)}

$$ B_5 := \sup_{I,J} \bigg(\iint\limits_{I\times J} v(x,y) dxdy \bigg) |I|^{q(\alpha-1/p)}
|J|^{q(\beta-1/p)} <\infty,$$ where $I$ and $J$ are arbitrary
bounded intervals in ${\Bbb{R}}$.} \vskip+0.2cm

 Exploring
the two-weight problem for the strong fractional maximal function
of variable order, in particular, we prove an analog of Theorem D
in variable exponent Lebesgue spaces  when the exponent of the
initial Lebesgue space is constant.

Let $p$ be a non--negative measurable function on ${\Bbb{R}}^n$.
Suppose that $E$ is a measurable subset of ${\Bbb{R}}^n$. In the
sequel we will use the following notation:
\begin{eqnarray*}
 p_-(E):= \inf_{E} p; \;\; p_+(E):= \sup_{E} p; \;\;\  p_-:= p_-({\Bbb{R}}^n); \;\; p_+:=
 p_+({\Bbb{R}}^n).
\end{eqnarray*}

Let $\Omega$ be an open set in ${\Bbb{R}}$ and let the following condition holds for $p: \Omega \to
{\Bbb{R}}$:
$$ 1\leq p_-(\Omega) \leq p(t)\leq p_+(\Omega)<\infty, \;\; t\in {\Bbb{R}}^n. $$

By $L^{p(\cdot)}(\Omega)$ we denote the Banach space of
measurable functions $f: \Omega \to \mathbb{R}$ such that
$$
\|f\|_{L^{p(\cdot)}(\Omega)}:=\|f\|_{p(\cdot)}=\inf\left\{\lb>0:
\intl_\Omega \left|\frac{f(x)}{\lb}\right|^{p(x)} \;dx\le
1\right\}<\infty .
$$
In the sequel by $L_{w}^{p(\cdot)}(\Omega)$ is denoted the
weighted variable exponent Lebesgue space defined by the norm
$$
\|f\|_{L_{w}^{p(\cdot)}}:=\|fw\|_{L^{p(\cdot)}}.
$$

It is known (see e.g. \cite{KoRa}, \cite{Sa1}, \cite{KoSa1},
\cite{HaHaPe}) that $L^{p(\cdot)}$  is a Banach space. For other
properties of $L^{p(\cdot)}$ spaces  we refer e.g., to \cite{Sh},
\cite{KoRa}, \cite{Sa1}.

Weighted estimates for classical integral operators in variable
exponent Lebesgue spaces were investigated in the papers
\cite{KoSa1}-\cite{KoSa5}, \cite{SaVa}, \cite{SaSaVa}, \cite{DiSa}
for power--type weights and in \cite{EdKoMeJFSA}-\cite{EdKoMe2},
\cite{KoMe4}, \cite{KoSa3}, \cite{Kop},  \cite{MaHa}, \cite{MaZe}, \cite{DiHa}
for general--type weights. Moreover, in the latter paper a
complete solution of the one--weight problem for maximal functions
defined on Euclidean spaces are given in terms of
Muckenhoupt--type conditions.

Finally we point out that constants (often different constants in
the same series of inequalities) will generally be denoted by $c$
or $C$. Throughout the paper by the symbol $p'(x)$ is denoted the
function $p(x)/ (p(x)-1)$.

\vskip+0.1cm

\section{STRONG FRACTIONAL MAXIMAL FUNCTIONS IN $L^{p(\cdot)}$ SPACES}
\setcounter{equation}{0}

Let
$$
\left(M_{\alpha(\cdot),\beta(\cdot)}^{S}f\right)(x,y)=\sup_{\substack{Q\ni
x\\ J\ni y }}
|Q|^{\frac{\alpha(x)}{n}-1}|J|^{\frac{\beta(y)}{n}-1}\int\limits_{Q\times
J}|f(t, \tau)| dt d\tau, \;\;\; (x,y)\in {\Bbb{R}}^n \times {\Bbb{R}}^m,
$$
where $\alpha$ and $\beta$ are measurable functions on
$\mathbb{R}^{n}$ and $\mathbb{R}^{m}$ respectively satisfying the
conditions: $0<\alpha_{-}\leq\alpha_{+}<n$,
$0<\beta_{-}\leq\beta_{+}<n$ and the supremum is taken over all
cubes $Q\ni x$ and $J\ni y$ respectively in $\mathbb{R}^{n}$ and
$\mathbb{R}^{m}$.

If $\alpha(\cdot)\equiv\beta(\cdot)$, then we denote
$M^{S}_{\alpha(\cdot),\beta(\cdot)}$ by $M^{S}_{\alpha(\cdot)}$.
Further,  if $\alpha\equiv const$, then we use the symbol
$M^{S}_{\alpha}$ instead of $M_{\alpha(\cdot)}^{S}$.

Let us take the case $n=m=1$ and consider the operator
$M^{S}_{\alpha}$ with constant parameter $\alpha$:
$$
\left(M_{\alpha}^{S}f\right)(x,y)=\sup_{\substack{I\ni x \\ J\ni
y}} |I\times J|^{\alpha-1}\int\limits_{I\times
J}|f|,\;\;\;(x,y)\in\mathbb{R}^{2},
$$
where the supremum is taken over all bounded intervals $I\ni x$ and $J\ni
y$, $0<\alpha<1$.

\vskip+0.2cm

{\bf Theorem.}  {\em Let $1<p_{-}\leq p_{+}<\infty$ and let $0\leq\alpha<\frac{1}{p_{-}}$.
We set $q(x)=\frac{p(x)}{1-\alpha\cdot p(x)}$. Then
$M^{S}_{\alpha}$ is bounded from $L^{p(\cdot)}(\mathbb{R}^{2})$ to
$L^{q(\cdot)}(\mathbb{R}^{2})$ if and only if $p\equiv const$. }

\begin{proof} {\em Sufficiency} is obvious using iterating process of
one-dimensional $L^{p}\rightarrow L^{q}$ boundedness of the one--dimensional fractional maximal operator

$$
(M_{\alpha}f)(x)=\sup_{\substack{I\ni x \\ I\subset \mathbb{R} }}
\frac{1}{|I|^{1-\alpha}}\int\limits_{I}|f(t)| dt,\;\;\;0\leq\alpha<1.
$$

{\em Necessity.} We follow T. Kopaliani \cite{Kop1} which proved
the theorem for $\alpha=0$. First we see that if $M^{S}_{\alpha}$
is bounded from $L^{p(\cdot)}(\mathbb{R}^{2})$ to
$L^{q(\cdot)}(\mathbb{R}^{2})$, then
$$
\sup\limits_{R}A_{R}:=\sup\limits_{R}\frac{1}{|R|^{1-\alpha}}\|\chi_{R}\|_{L^{q(\cdot)}}\|\chi_{R}\|_{L^{p'(\cdot)}}<\infty,
$$
where the supremum is taken over all rectangles $R$ in
$\mathbb{R}^{2}$.

Indeed, let $\|f\|_{L^{p(\cdot)}(\mathbb{R}^{2})}\leq 1$. Then for
every rectangle $R$ we have

$$
c\geq \|M^{S}_{\alpha}f\|_{L^{q(\cdot)}(\mathbb{R}^{2})}\geq
\|M^{S}_{\alpha}f\|_{L^{q(\cdot)}(R)} \geq
\|\chi_{R}\|_{L^{q(\cdot)}}|R|^{\alpha-1}\iint\limits_{R}|f(t, \tau)|dt d\tau.
$$
Taking now the supremum with respect to $f$,
$\|f\|_{L^{p(\cdot)}}\leq 1$, we find that
$$
|R|^{\alpha-1}\|\chi_{R}\|_{L^{q(\cdot)}}\|\chi_{R}\|_{L^{p'(\cdot)}}\leq
c
$$
for all $R\subset\mathbb{R}^{2}$.

Further, suppose the contrary: $p$ is not constant, i. e.
$$
\inf\limits_{\mathbb{R}^{2}}p(t)<\sup\limits_{\mathbb{R}^{2}}p(t).
$$
By Luzin's theorem, there is a family of pointwise disjoint sets
$F_{i}$ satisfying the conditions:

$(i)$ $|\mathbb{R}^{2}\backslash \cup_{j}F_{j}|=0$;

$(ii)$ functions $p:F_{i}\rightarrow\mathbb{R}$ are continuous;

$(iii)$ for every fixed $i$, all points of $F_{i}$ are points of
density with respect to the basis consisting of all open
rectangles in $\mathbb{R}^{2}$.

We can find a pair of the type $\left( (x_{0},y_{1}),(x_{0},y_{2})
\right)$ or $\left( (x_{1},y_{0}),(x_{2},y_{0}) \right)$ from
$\cup F_{i}$ such that $p(x_{0},y_{1})\neq p(x_{0},y_{2})$ or
$p(x_{1},y_{0})\neq p(x_{2},y_{0})$. Without loss of generality,
assume that this pair is $\left( (x_{0},y_{1}),(x_{0},y_{2})
\right)$ such that $(x_{0},y_{1})\in F_{1}$ and $(x_{0},y_{2})\in
F_{2}$, $y_{1}<y_{2}$.

Let $0<\varepsilon<1$ be fixed number. Then there is a number
$\delta>0$ such that for any rectangles $Q_{1}\ni (x_{0},y_{1})$
and $Q_{2}\ni (x_{0},y_{2})$ with diameters less than $\delta$,
the following inequalities hold:

$$
|Q_{1}\cap F_{1}|>(1-\varepsilon)|Q_{1}|,\;\;\;|Q_{2}\cap
F_{2}|>(1-\varepsilon)|Q_{2}|, \eqno{(1.1)}
$$

$$
p_{Q_{1}}=\sup\limits_{Q_{1}\cap
F_{1}}p(x,y)<c_{1}<c_{2}<\inf\limits_{Q_{2}\cap
F_{2}}p(x,y)=p_{Q_{2}}, \eqno{(1.2)}
$$

where $c_{1}$ and $c_{2}$ are some constants.

Let $Q_{1,\tau}$ and $Q_{2,\tau}$ be rectangles with properties
(1.1) and (1.2) with the forms
$(x_{0}-\tau,x_{0}+\tau)\times(a,b)$ and
$(x_{0}-\tau,x_{0}+\tau)\times(c,d)$ respectively, where
$a<b<c<d$.

Observe now that the following embeddings hold:
$$
L^{q(\cdot)}(Q_{2,\tau})\hookrightarrow L^{q_{Q_{2}}}(Q_{2,\tau})
$$
$$
L^{p'(\cdot)}(Q_{1,\tau})\hookrightarrow
L^{(p_{Q_{1}})'}(Q_{1,\tau}),
$$
where $q_{Q_{2}}=\inf\limits_{Q_{2}}q=\frac{p_{Q_{2}}}{1-\alpha
p_{Q_{2}}}$, $(p_{Q_{1}})'=\frac{p_{Q_{1}}}{p_{Q_{1}}-1}$.
Further, for the rectangle $Q_{\tau}=(x_{0}-\tau,x_{0}+\tau)\times
(a,d)$, we have that
$$
A_{\tau}:=\frac{1}{|Q_{\tau}|^{1-\alpha}}
\|\chi_{Q_{\tau}}\|_{L^{q(\cdot)}}\|\chi_{Q_{\tau}}\|_{L^{p'(\cdot)}}
$$
$$
\geq\frac{1}{[2\tau(d-a)]^{1-\alpha}}\|\chi_{Q_{2,\tau\cap
F_{2}}}\|_{L^{q(\cdot)}} \|\chi_{Q_{1,\tau\cap
F_{1}}}\|_{L^{p'(\cdot)}}
$$
$$
\geq\frac{c}{[2\tau(d-a)]^{1-\alpha}}[2\tau(d-c)]^{\frac{1}{q_{Q_{2}}}}
[2\tau(b-a)]^{1-\frac{1}{p_{Q_{1}}}}
$$
$ =c \tau^{\alpha-1+\frac{1}{q_{Q_{2}}}+1-\frac{1}{p_{Q_{1}}}}=c
\tau^{\alpha-\left[\frac{1}{p_{Q_{1}}}+\frac{1}{q_{Q_{2}}}\right]}\rightarrow
\infty $ as $\tau\rightarrow 0$ because
$\alpha-\frac{1}{p_{Q_{1}}}+\frac{1}{q_{Q_{2}}}=\alpha-\frac{1}{p_{Q_{1}}}+\frac{1}{p_{Q_{2}}}-\alpha<0$
(recall that $a$, $b$, $c$ and $d$ are fixed).

This contradicts the condition
$$
\sup_{R} A_{R}<\infty.
$$
\end{proof}

\section{Double Hardy Transform in $L_{w}^{p(\cdot)}$ spaces}

Let
$$
(\mathcal{H}_{2}f)(x,y)=\int\limits_{0}^{x}\int\limits_{0}^{y}f(t,\tau)dt
d\tau,\;\;\;(x,y)\in\mathbb{R}^{2}_{+}.
$$


First we prove the following lemma: \vskip+0.1cm
{\bf Lemma 3.1.} {\em  Let $p$ be a constant satisfying the
condition $1<p<\infty$. Suppose that $0<b\leq \infty$. Let $\rho$
be an almost everywhere positive function on $[0,b)$. Then there
is a positive constant $c$ such that for all $f\in
L^p_{\rho}([0,b))$, $f\geq 0$, the inequality
$$  \int\limits_0^b \bigg( \frac{1}{\lambda([0,x]}
\int\limits_0^x f(t) dt \bigg)^{p} \lambda(x) dx \leq C
\int\limits_0^b (f(x)\rho(x))^p dx $$ holds, where $\lambda(x)=
\rho^{-p'}(x)$ and $\lambda([0,x]):=\int\limits_{0}^x\lambda(t)
dt$.}

\vskip+0.2cm

\begin{proof} It is enough to show that (see e.g. \cite{Maz}, Chapter 1)
the condition
$$ \sup\limits_{0<t<b}\bigg( \int_t^b \lambda([0,x])^{-p}
\lambda(x) dx\bigg)\bigg(\int_0^t \lambda(x) dx
\bigg)^{p-1}<\infty $$ is satisfied.

To check that this condition holds observe that
$$ \int\limits_t^b \lambda([0,x])^{-p} \lambda(x) dx =
\int\limits_t^b \bigg( \int\limits_0^x\lambda(\tau) d\tau
\bigg)^{-p} \lambda(x) dx $$

$$ = \frac{1}{1-p} \int\limits_t^b d \bigg( \int\limits_0^x
\lambda(\tau) d\tau\bigg)^{1-p} = \frac{1}{p-1} \Bigg[ \bigg(
\int\limits_0^t \lambda(\tau) d\tau \bigg)^{1-p} - \bigg(
\int\limits_0^b \lambda(\tau) d\tau \bigg)^{1-p}\Bigg] $$

$$\leq \frac{1}{p-1} \Bigg[ \bigg( \int\limits_0^t \lambda(\tau) d\tau
\bigg)^{1-p} + \bigg( \int\limits_0^b \lambda(\tau) d\tau
\bigg)^{1-p}\Bigg].$$ Now the result follows easily.

\end{proof}
\vskip+0.2cm

{\bf Theorem 3.1.} {\em Let $p$ be constant and let $1<p\leq q_{-}\leq q_{+}<\infty$. Suppose that $v$ and $w$ be weights
on $\mathbb{R}^{2}_{+}$ with
$w(x,y)=w_{1}(x)w_{2}(y)$ for some one--dimensional weights $w_1$
and $w_2$. Then $\mathcal{H}_{2}$ is bounded from
$L_{w}^{p}(\mathbb{R}^{2}_{+})$ to
$L_{v}^{q(\cdot)}(\mathbb{R}^{2}_{+})$ if and only if
$$
B:=\sup\limits_{a,b>0}\|v(\chi_{J_{ab}^{\infty}})\|_{L^{q(\cdot)}(\mathbb{R}^{2}_{+})}
\bigg\|\frac{1}{w}\cdot\chi_{J^{0}_{ab}}\bigg\|_{L^{p'}(\mathbb{R}^{2}_{+})}<\infty,
$$
where $J_{ab}^{\infty}=[a,\infty)\times[b,\infty)$ and
$J_{ab}^{0}=[0,a)\times [0,b)$. }

\begin{proof} Necessity follows by the standard way taking the
test function
$$
f(x,y)=\left(\int\limits_{0}^{a}\int\limits_{0}^{b}w^{-p'}(x,y)dxdy\right)\chi_{[0,a]\times
[0,b]}(x,y)
$$
in the two-weight inequality and do simple estimates.

{\em Sufficiency}. Suppose that $f\geq 0$ and $\|f
\|_{L^{p}_w({\Bbb{R}}^2)}\leq 1$. Let $\{x_{k}\}$ and $\{y_{j}\}$
be sequences of positive numbers chosen so that
$$
\int\limits_{0}^{x_{k}}w_{1}^{-p'}=2^{k},\;\;\;\;\;\int\limits_{0}^{y_{j}}w_{2}^{-p'}=2^{j}.
\eqno{(3.1)}
$$
Without loss of generality assume that
$\int\limits_{0}^{\infty}w_{1}^{-p'}=\int\limits_{0}^{\infty}w_{2}^{-p'}=\infty$.
Then
$[0,\infty)=\bigcup\limits_{k}[x_{k},x_{k+1})=\bigcup\limits_{j}[y_{j},y_{j+1})$.
On the other hand,
$\mathbb{R}^{2}_{+}=\bigcup\limits_{k,j}[x_{k},x_{k+1})\times
[y_{j},y_{j+1})$. It is easy to see that equalities (3.1) imply:
$$
\int\limits_{x_{k}}^{x_{k+1}}w_{1}^{-p'}=2^{k},\;\;\;\;\;\int\limits_{y_{j}}^{y_{j+1}}w_{2}^{-p'}=2^{j}.
\eqno{(3.2)}
$$
Denote: $[x_{k},x_{k+1})=:E_{k},$, $[y_{j},y_{j+1})=:F_{j}$.

Let us choose $r$ so that $p\leq r\leq q_{-}$. Then
$$
\|v(\mathcal{H}_{2}f)\|^{r}_{L^{q(\cdot)}(\mathbb{R}^{2}_{+})}=\|[v(\mathcal{H}_{2}f)]^{r}\|_{L_{v}^{q(\cdot)/r}
(\mathbb{R}^{2}_{+})}
$$
$$
\leq c\sup\limits_{\|h\|_{L^{(q(\cdot)/r)'}\leq 1}}
\iint_{\mathbb{R}^{2}}(v(x,y))^{r}(\mathcal{H}_{2}f(x,y))^{r}h(x,y)dxdy.
$$
Further, taking (3.1) and (3.2) into account we have that
$$
\iint_{\mathbb{R}^{2}_{+}}(v(x,y))^{r}(\mathcal{H}_{2}f)^{r}(x,y)h(x,y)dxdy
$$
$$
=\sum\limits_{k,j}\left[\int\limits_{x_{k}}^{x_{k+1}}\int\limits_{y_{j}}^{y_{j+1}}v^{r}(x,y)h(x,y)dxdy\right]
\left[\int\limits_{0}^{x_{k+1}}\int\limits_{0}^{y_{j+1}}f\right]^{r}
$$
$$
\leq\sum\limits_{k,j}\|v^{r}(\cdot)\|_{L^{q(\cdot)/r} (E_{k}\times
F_{j})}\|h\|_{L^{(q(\cdot)/r)'}(\mathbb{R}^{2}_{+})}
\left[\int\limits_{0}^{x_{k+1}}\int\limits_{0}^{y_{j+1}}f\right]^{r}
$$
$$
\leq\sum\limits_{k,j}\|v\|^{r}_{L^{q(\cdot)}(E_{k}\times
F_{j})}\left[\int\limits_{0}^{x_{k+1}}\int\limits_{0}^{y_{j+1}}f\right]^{r}
$$
$$
\leq
B^{r}\sum\limits_{k,j}\|w_{1}^{-1}\|^{-r}_{L^{p'}([0,x_{k}])}\|w_{2}^{-1}\|^{-r}_{L^{p'}([0,y_{j}])}
\left[\int\limits_{0}^{x_{k+1}}\int\limits_{0}^{y_{j+1}}f\right]^{r}
$$
$$
=c_{r}B^{r}\sum\limits_{k,j}\left[\int\limits_{x_{k+1}}^{x_{k+2}}
\int\limits_{y_{j+1}}^{y_{j+2}}(w_{1}(x)w_{2}(y))^{-p'}dxdy\right]^{\frac{r}{p}}
\cdot\left[\frac{1}{\sigma_{1}(E_{k})\sigma_{2}(F_{j})}\int\limits_{0}^{x_{k+1}}\int\limits_{0}^{y_{j+1}}f\right]^{r}
$$
(where $\sigma_{1}(E_{k}):=\int\limits_{E_{k}}w_{1}^{-p'}$,
$\sigma_{2}(F_{j}):=\int\limits_{F_{j}}w_{2}^{-p'}$)
$$
\leq c
B^{r}\sum\limits_{k,j}\left(\int\limits_{x_{k+1}}^{x_{k+2}}\int\limits_{y_{j+1}}^{y_{j+2}}
(w_{1}(x)w_{2}(y))^{-p'}\left[\frac{1}{\sigma_{1}([0,x_{k+2}])\sigma_{2}([0,y_{j+2}])}\int\limits_{0}^{x_{k+1}}
\int\limits_{0}^{y_{j+1}}f\right]^{p}dxdy\right)^{r/p}
$$
$$
\leq c
B^{r}\sum\limits_{k,j}\left[\int\limits_{x_{k+1}}^{x_{k+2}}\int\limits_{y_{j+1}}^{y_{j+2}}
[w_{1}(x)w_{2}(y)]^{-p'}\left[\frac{1}{\sigma_{1}([0,x])\sigma_{2}([0,y])}\int\limits_{0}^{x}\int\limits_{0}^{y}
f\right]^{p}dxdy\right]^{r/p}
$$
$$
\leq c
B^{r}\left[\iint\limits_{\mathbb{R}^{2}_{+}}[w_{1}(x)w_{2}(y)]^{-p'}\left[\frac{1}{\sigma_{1}([0,x])\sigma_{2}[0,y]}
\int\limits_{0}^{x}\int\limits_{0}^{y}f\right]^{p}dxdy\right]^{r/p}.
$$
Observe now that Lemma 3.1 implies the inequality
$$
\int\limits_{\mathbb{R}_{+}}\left[\frac{1}{\sigma([0,x])}\int\limits_{0}^{x}f\right]^{p}d\sigma(x)\leq
c\int\limits_{\mathbb{R}_{+}}(f(x)w(x))^{p}dx. \eqno{(3.3)}
$$
By using inequality (3.3) twice together with Fubini's theorem we
find that
$$
\left[\iint\limits_{\mathbb{R}_{+}^{2}}w^{-p'}(x,y)\left[\frac{1}{\sigma([0,x]\times
[0,y])}\int\limits_{0}^{x}\int\limits_{0}^{y}f\right]^{p}dxdy\right]^{r/p}
$$
$$
\leq
c\left[\iint\limits_{\mathbb{R}_{+}^{2}}[f(x,y)]^{p}(w(x,y))^{p}dxdy\right]^{\frac{r}{p}}\leq
c.
$$
\end{proof}

{\bf Corollary 3.1.} (Trace inequality) {\em  Let $1<p\leq
q_{-}\leq q_{+}<\infty$ and let $v$ be a. e. positive function on
$\mathbb{R}^{2}_{+}$. Then $\mathcal{H}_{2}$ is bounded from
$L^{p}(\mathbb{R}^{2}_{+})$ to $L^{q(\cdot)}_{v}(\mathbb{R}^{2}_+)$
if and only if
$$
\sup\limits_{a,b>0}\|v\chi_{J_{ab}^{\infty}}\|_{L^{q(\cdot)}(\mathbb{R}^{2}_{+})}
(ab)^{\frac{1}{p'}}<\infty.
$$}
\vskip+0.1cm

{\bf Definition 3.1.}  Let $\Omega$ be an open set in ${\Bbb{R}}^n$. We say that the exponent function
$p(\cdot)\in\mathcal{P}(\Omega)$ if there is a constant
$0<\delta<1$ such that
$$
\int\limits_{\Omega}\delta^{\frac{p(x)p_{-}}{p(x)-p_{-}}}dx<+\infty.
$$
Further, we say that $p(\cdot)\in\mathcal{P}_{\infty}(\Omega)$ if
$$
|p(x)-p(y)|\leq\frac{c}{\ln(e+|x|)}
$$
for all $x,y\in\Omega$ with $|y|\geq |x|$.

\vskip+0.1cm

{\bf Corollary 3.2.} {\em  Let $1<p_{-}\leq q_{-}\leq
q_{+}<\infty$ with $p_{+}<\infty$. Let $v$ and $w$ be a. e.
positive functions on $\mathbb{R}^{2}$ with
$w(x,y)=w_{1}(x)w_{2}(y)$. Suppose that
$p\in\mathcal{P}(\mathbb{R}^{2}_{+})$. If
$$
\sup\limits_{a,b>0}\big\|v\chi_{J_{ab}^{\infty}}\big\|_{L^{q(\cdot)}(\mathbb{R}^{2}_{+})}
\bigg\|\frac{1}{w}\cdot\chi_{J_{ab}^{0}}\bigg\|_{L^{(p_{-})'}(\mathbb{R}^{2}_{+})}<\infty,
\eqno{(3.4)}
$$
then $\mathcal{H}_{2}$ is bounded from
$L_{w}^{p(\cdot)}(\mathbb{R}^{2}_+)$ to
$L_{v}^{q(\cdot)}(\mathbb{R}^{2}_+)$. }

\begin{proof}
Recall  that (see \cite{Di1}) if
$p\in\mathcal{P}(\mathbb{R}^{2}_{+})$, then
$L^{p(\cdot)}(\mathbb{R}_{+}^{2})\hookrightarrow
L^{p_{-}}(\mathbb{R}^{2}_{+})$. Now Theorem 3.1 completes the
proof.
\end{proof}

\vskip+0.2cm
{\bf Corollary 3.3.} {\em  Let $1<p_{-}\leq q_{-}\leq
q_{+}<\infty$ with $p_{+}<\infty$ and $p_{\infty}=p_{-}$. Assume
that $p\in\mathcal{P}_{\infty}(\mathbb{R}^{2}_{+})$. Suppose that
$v$ and $w$ are a. e. positive functions on $\mathbb{R}^{2}_{+}$
with $w(x,y)=w_{2}(x)w_{2}(y)$. If $(3.4)$ holds, then
$\mathcal{H}_{2}$ is bounded from
$L_{w}^{p(\cdot)}(\mathbb{R}^{2})$ to
$L_{v}^{q(\cdot)}(\mathbb{R}^{2})$. } \vskip+0.2cm

 Let us now consider the
operator $\mathcal{H}_{2}$ on a rectangle $J:=[0,a_{0}]\times
[0,b_{0}]$. In the sequel the following notation will be used:
$$
J_{ab}^{0}:=[0,a]\times [0,b],\;\;\;J_{ab}^{1}:=[a,a_{0}]\times
[b,b_{0}].
$$
The arguments used in the proof of Theorem 3.1 enable us to
formulate the next statements: \vskip+0.2cm

{\bf Theorem 3.2.} {\em  Let $1<p_{-}(J)\leq q_{-}(J)\leq
q_{+}(J)<\infty$ with $p_{+}(J)<\infty$. Suppose that $v$ and $w$
are a. e. positive functions on $J$ with $w(x,y)=w_{1}(x)w_{2}(y)$
for some one-dimensional weights $w_1$ and $w_2$. If
$$
\sup_{\substack{0<a\leq a_{0} \\ 0<b\leq b_{0}
}}\|v\chi_{J_{ab}^{1}}\|_{L^{q(\cdot)}(\mathbb{R}^{2}_{+})}
\|w^{-1}\chi_{J_{ab}^{0}}\|_{L^{(p_{-})'}(\mathbb{R}^{2}_{+})}<\infty
$$
then $\mathcal{H}_{2}$ is bounded from $L_{w}^{p(\cdot)}(J)$ to
$L_{v}^{q(\cdot)}(J).$ }

\vskip+0.1cm

{\bf Corollary 3.4.} {\em  Let $1<p_{-}(J)\leq q_{-}(J)\leq
q_{+}(J)<\infty$ with $p_{-}(J)=p(0)$ and $p_{+}<\infty$. Let $v$
be a. e. positive function on $J$. Then $\mathcal{H}_{2}$ is
bounded from $L^{p(\cdot)}(J)$ to $L^{q(\cdot)}_{v}(J)$ if
$$
\sup_{\substack{0<a\leq a_{0} \\ 0<b\leq b_{0} }}
\|v\chi_{J_{ab}^{1}}\|_{L^{q(\cdot)}(\mathbb{R}^{2}_{+})}
\|\chi_{J_{ab}^{0}}\|_{L^{p'(0)}(\mathbb{R}^{2}_{+})}<\infty
$$}

\vskip+0.1cm

{\bf Corollary 3.5.} {\em  There is non-constant exponent $p$ on
$[0,2]^2$ such that the double average operator
$$
(Af)(x,y)=\frac{1}{xy}\int\limits_{0}^{x}\int\limits_{0}^{y}
f(t,\tau)dtd\tau
$$
is bounded in $L^{p(\cdot)}([0,2]^{2})$. }

\begin{proof}
Let $p$ be defined as follows:
$$
p(x,y)=\left\{%
\begin{array}{ll}
    3, & \hbox{$(x,y)\in [1,2]^{2}$;} \\
    2, & \hbox{$(x,y)\in [0,2]^{2}\backslash [1,2]^{2}$} .\\
\end{array}%
\right.
$$
It is clear that $p(0,0)=p_{-}=2$.

Also, it is easy to check that
$$
\sup\limits_{0<a,b\leq 2} \|(xy)^{-1}\chi_{[a,2]\times
[b,2]}(x,y)\|_{L^{p(x,y)}(\mathbb{R}^{2}_{+})}(ab)^{\frac{1}{p'(0,0)}}<\infty.
$$
Corollary 3.4 completes the proof.
\end{proof}

\section{TWO-WEIGHT ESTIMATES FOR STRONG FRACTIONAL MAXIMAL FUNCTION IN $L^{p(\cdot)}$ SPACES}

\

In order to establish two-weight estimates for strong fractional
maximal function of variable order we need the next
Carleson-H\"{o}rmander's embedding theorem regarding dyadic
intervals.

A weight function $\rho$ is said to be satisfying the dyadic
reverse doubling condition $(\rho\in RD^{(d)}(\mathbb{R}))$ if for
any two dyadic intervals $I$ and $I'$ with $I\subset I'$,
$|I|=\frac{|I'|}{2}$ the inequality
$$
\rho(I')\leq b\rho(I)
$$
holds with some constant $b>1$.
\vskip+0.1cm

{\bf Theorem E (\cite{Ta}, \cite{SaWh}, Lemma 3.10).}  {\em Let  $p$ and $q$ be constants satisfying the condition $1<p<q<\infty$ and let $\rho$ be a
weight function on $\mathbb{R}$ such that $\rho^{1-p'}$ satisfies
the dyadic reverse doubling condition. Let $\{c_{I}\}$ be a sequence of
non--negative numbers corresponding to dyadic intervals $I$ in
$\mathbb{R}$. Then the following two statements are equivalent:

$\rm{(i)}$\;There is a positive constant $C$ such that
$$
\sum\limits_{I\in
\mathcal{D}}c_{I}\left(\frac{1}{|I|}\int\limits_{I}g(x)dx\right)^{q}
\leq c\left(\int\limits_{\mathbb{R}}g(x)^{p}\rho(x)dx\right)^{q/p}
$$
for all nonnegative $g\in L_{\rho}^{p}(\mathbb{R})$;

$\rm{(ii)}$\;There is a positive constant $C_{1}$ such that
$$
c_{I}\leq
C_{1}|I|^{q}\left(\int\limits_{I}\rho(x)^{1-p'}dx\right)^{-q/p'}
$$
for all $I\in D$.}
\vskip+0.2cm

This result yields the following corollary.
\vskip+0.1cm
{\bf Corollary A.} {\em Let $p$ and $q$ be constants satisfying the condition $1<p<q<\infty$ and let $\rho$ be a weight
function on $\mathbb{R}$ such that $\rho^{1-p'}$ satisfies the
dyadic reverse doubling condition. Then the Carleson-H\"{o}rmander
inequality
$$
\sum\limits_{I\in
\mathcal{D}}\left(\int\limits_{I}\rho^{1-p'}(x)dx\right)^{-q/p'}
\left(\int\limits_{I}f(x)dx\right)^{q} \leq
c\left(\int\limits_{\mathbb{R}}f^{p}(x)\rho(x)dx\right)^{q/p}
$$
holds for all nonnegative $f\in L_{\rho}^{p}(\mathbb{R})$.}

Recall that (see Section 2)
$$
\Big(M_{\al(\cdot),\bt(\cdot)}^S
f\Big)(x,y)=\sup_{\begin{subarray}{c}
Q\in x\\
J\ni y
\end{subarray}}
|Q|^{\frac{\al(x)}{n}-1}|J|^{\frac{\bt(y)}{m}-1}\iint_{Q\times
J}|f|,
$$
where $(x,y)\in\bR^n\times\bR^m$ and $Q$ and $J$ are cubes in
$\bR^n$ and $\bR^m$ respectively. For simplicity we take $n=m=1$
and consider the strong maximal operator
$$
\Big(M_{\al(\cdot),\bt(\cdot)}^S
\Big)f(x,y)=\sup_{\begin{subarray}{c}
I\ni x\\
J\ni y
\end{subarray}}
|I|^{\al(x)-1}|J|^{\bt(y)-1}\iint_{I\times J}|f|,\;(x,y)\in \bR^2,
$$
where $0<\al_-\leq \al_+<1$, $0<\bt_-\leq\bt_+<1$ and $I$ and $J$
are intervals in $\bR$.

Together with the operator $M_{\al(\cdot),\bt(\cdot)}^S$ we are
interested in the dyadic strong fractional maximal operator
$$
\Big(M_{\al(\cdot),\bt(\cdot)}^{S,(d)}f
\Big)(x,y)=\sup_{\begin{subarray}{c}
I\ni x\\
J\ni y\\
I,J\in D(\bR)
\end{subarray}}
|I|^{\al(x)-1}|J|^{\bt(y)-1}\iint_{I\times J}|f(t, \tau)|dt d\tau,\;\;(x,y)\in \bR^2,
$$
where $I$ and $J$ belong to the dyadic lattice $D(\bR)$ of $\bR$.

\subsection*{The Fefferman-Stein Type Inequalities. Criteria for the  Trace
Inequality}

We start by the Fefferman-Stein type inequality. The original
inequality for fractional maximal operator defined on cubes in
$L^p$ spaces with constant $p$ was derived by E. T. Sawyer.

\vskip+0.2cm

{\bf Theorem 4.1.} {\em  Let $1<p_-\leq p_+<q_-\leq q_+<\infty$
and let
$\frac{1}{p_-}-\frac{1}{q_+}<\al_-\leq\al_+<\frac{1}{p_-}$, $\frac{1}{p_-} -\frac{1}{q_+}< \bt_- \leq \bt_+<\frac{1}{p_-}$.
Then there is a positive constant $c$ such that
$$
\|(M_{\al(\cdot),\bt(\cdot)}^Sf)v\|_{L^{q(\cdot)}(\bR^2)}\leq
c\|f(\wt{M}_{\al(\cdot),\bt(\cdot)}v)\|_{L^{p(\cdot)}(\bR^2)},
$$
where \begin{gather*}
(\widetilde{M}_{\al(\cdot),\bt(\cdot)}v)(x,y)=\max\{(\wt{M}_{\al(\cdot),\bt(\cdot)}^{(1)}v)(x,y)
(\wt{M}_{\al(\cdot),\bt(\cdot)}^{(2)}v)(x,y)\},\\
(\wt{M}_{\al(\cdot),\bt(\cdot)}^{(1)}v)(x,y)=
\sup_{\begin{subarray}{c}
I\ni x\\
J\ni y
\end{subarray}}
|I\times J|^{-\frac{1}{p_-}}\|v(\cdot)|I|^{\al(\cdot)}|J|^{\bt(\cdot)}\|_{L^{q(\cdot)}(I\times J)},\\
(\wt{M}_{\al(\cdot),\bt(\cdot)}^{(2)}v)(x,y)=
\sup_{\begin{subarray}{c}
I\ni x\\
J\ni y
\end{subarray}}
|I\times
J|^{-\frac{1}{p_+}}\|v(\cdot)|I|^{\al(\cdot)}|J|^{\bt(\cdot)}\|_{L^{q(\cdot)}(I\times
J)},
\end{gather*}
} \vskip+0.2cm

{\bf Corollary 4.1.} {\em  Let $p$ be constant such that
$1<p<q_-\leq q_+<\infty$ and let
$\frac{1}{p}-\frac{1}{q_+}<\al_-\leq\al_+<\frac{1}{p}$,
$\frac{1}{p}-\frac{1}{q_-}<\bt_-\leq\bt_+<\frac{1}{p}$. Then the
following inequality holds:
$$
\|(M_{\al(\cdot),\bt(\cdot)}^Sf)v\|_{L^{q(\cdot)}(\bR^2)}\leq
c\|f(\wt{M}_{\al(\cdot),\bt(\cdot)}v)\|_{L^{p}(\bR^2)},
$$
where
$$(\wt{M}_{\al(\cdot),\bt(\cdot)}v)(x,y)=
\sup_{\begin{subarray}{c}
I\ni x\\
J\ni y
\end{subarray}}
|I\times
J|^{-\frac{1}{p}}\|v(\cdot)|I|^{\al(\cdot)}|J|^{\bt(\cdot)}\|_{L^{q(\cdot)}(I\times
J)}.
$$
}

\vskip+0.1cm

{\em Remark} 4.1.  Notice that for $\al\equiv\const$, $\bt\equiv
\const$, the operator $\wt{M}_{\al,\bt}$ has the form
$$
(\wt{M}_{\al,\bt}v)(x,y)= \sup_{\begin{subarray}{c}
I\ni x\\
J\ni y
\end{subarray}}
|I|^{\al-\frac{1}{p}}|J|^{\bt-\frac{1}{p}}\|v(\cdot)\|_{L^{q(\cdot)}(I\times
J)}.
$$
\vskip+0.1cm

{\em Remark} 4.1. If $q=\const$, then
$$
(\wt{M}_{\al(\cdot),\bt(\cdot)}v)(x,y)=\sup_{\begin{subarray}{c}
I\ni x\\
J\ni y
\end{subarray}}|I\times J|^{-\frac{1}{p}}\bigg(\iint_{I\times J}v^q(x,y)|I|^{q\al(x)}|J|^{q\bt(y)}dxdy\bigg)^{\frac{1}{q}}
$$
\vskip+0.2cm

{\bf Corollary 4.2.} [{Trace inequality}] {\em Let $1<p_-\leq
p_+<q_-\leq q_+<\infty$ and let
$\frac{1}{p_-}-\frac{1}{q_+}<\al_-\leq\al_+<\frac{1}{p_-}$,
$\frac{1}{p_-}-\frac{1}{q_+}<\bt_-\leq\bt_+<\frac{1}{p_-}$.
Suppose that the weight function $v$ satisfies the condition
$$
\sup_{I,J\sbs\bR}\||I|^{\al(\cdot)}|J|^{\bt(\cdot)}v(\cdot)\|_{L^{q(\cdot)}(I\times
J)}|I\times J|^{-\frac{1}{\overline{p}_{I\times J}}}<\infty,
$$
where
$$
\ol{p}_{I\times J}=\begin{cases}
p_-\;\;\;\text{if}\;\;\;|I||J|\leq 1\\
p_+\;\;\;\text{if}\;\;\;|I||J|>1
\end{cases}
$$
Then $M_{\al(\cdot),\bt(\cdot)}^S$ is bounded from
$L^{p(\cdot)}(\bR^2)$ to $L_v^{q(\cdot)}(\bR^2)$. }

\vskip+0.2cm

{\bf Theorem 4.2.} [{Criteria for the trace inequality}] {\em Let
$p$ be constant and let $1<p<q_-\leq q_+<\infty$. Suppose that
$\frac{1}{p}-\frac{1}{q_+}<\al_-\leq\al_+<\frac{1}{p}$,
$\frac{1}{p}-\frac{1}{q_+}<\bt_-\leq\bt_+<\frac{1}{p}$. Then
$M_{\al(\cdot),\bt(\cdot)}^S$ is bounded from $L^p(\bR^2)$ to
$L_v^{q(\cdot)}(\bR^2)$ if and only if
$$
\sup_{I,J\sbs\bR}\||I|^{\al(\cdot)}|J|^{\bt(\cdot)}v(\cdot)\|_{L^{q(\cdot)}(I\times
J)}|I\times J|^{-\frac{1}{p}}<\infty.
$$}
\vskip+0.2cm

{\bf Theorem 4.3.} {\em  Let $p$ be constant and let $1<p<q_-\leq
q_+<\infty$. Suppose that $0<\al_-\leq \al_+<1$, $0<\beta_-\leq
\beta_+<1$. Let $v$ and $w$ be weight functions in $\bR^2$ and let
$w$ is of product type, i.e. $w(x,y)=w_1(x)w_2(y)$. Then
$M_{a(\cdot),\bt(\cdot)}^s$ is bounded from $L^p_{w}(\bR^2)$ to
$L^{q(\cdot)}_{v}(\bR^2)$ if and only if
$$
\sup_{I,J\sbs \bR^2}(|I| |J|)^{-1}\|v(\cdot)|I|^{\al(\cdot)}
|J|^{\bt(\cdot)}\|_{L^{q(\cdot)}(I\times
J)}\|w^{-1}\|_{L^{p'(\cdot)}(I\times J)}<+\infty
$$
provided that $w_1^{-p'}, w_2^{-p'}\in RD^{(d)}(\bR)$. }
\vskip+0.2cm

{\bf Corollary 4.3.} {\em  Let $1<p_-<q_-\leq q_+<\infty$ with
$p_+<\infty$, $0<\al_-\leq\al_+<1$, $0<\bt_-\leq \bt_+<1$. Suppose
that $p\in \mathcal{P}(\bR^2)$. Assume that $v$ and $w$ are weight
functions on $\bR^2$ and that $w(x,y)=w_1(x)w_2(y)$ with
$w_1^{-(p_-)'}, w_2^{-(p_-)'}\in R D^{(d)}(\bR)$. If the condition
$$
\sup_{I,J\sbs \bR^2}(|I| |J|)^{-1}\|v(\cdot)|I|^{\al(\cdot)}
|J|^{\bt(\cdot)}\|_{L^{q(\cdot)}(I\times
J)}\|w^{-1}\|_{L^{(p_-)'}(I\times J)}<+\infty
\eqno{(4.1)} $$
holds, then $M_{\al(\cdot),\bt(\cdot)}^S$ is bounded from
$L_w^{p(\cdot)}(\bR^2)$ to $L_v^{q(\cdot)}(\bR^2)$. } \vskip+0.2cm

\vskip+0.2cm

{\bf Corollary 4.4.} {\em  Let $1<p_-<q_-\leq q_+<\infty$ with
$p_+<\infty$, $0<\al_-\leq \al_+<1$, $0<\bt_-\leq \bt_+<1$.
Suppose that $p_-=p(\infty)$ and that $p\in
\mathcal{P}_\infty(\bR^2)$. Suppose that $v$ and $w$ are weights
on $\bR^2$ and $w(x,y)=w_1(x)w_2(y)$ with $w_1^{-(p_-)'}$,
$w_2^{-(p_-)'}\in RD^{(d)}(\bR)$. Then condition $(4.1)$
guarantees the boundedness of $M_{\al(\cdot),\bt(\cdot)}^s$, from
$L_w^{p(\cdot)}(\bR^2)$ to $L_v^{q(\cdot)}(\bR^2)$. }

\subsection*{Proofs of the Results}

\noindent{\it Proof of Theorem $4.1$.}  Recall that by
$M_{\al(\cdot),\bt(\cdot)}^{S,(d)}$ we denote the dyadic
fractional maximal operator. Without loss of generality we can
assume that $f\geq 0$ and $f$ is bounded with compact support.

It is obvious that for $(x,y)\in\bR^2$ there are dyadic intervals
$I(x)\ni x$, $J(y)\ni y$ such that
$$
\frac{2}{|I(x)|^{1-\al(x)}|J(y)|^{1-\bt(y)}}\iint_{I(x)\times
J(y)}|f(t, \tau)| dt d\tau >(M_{\al(\cdot),\bt(\cdot)}^{S,(d)}f)(x,y).
$$

Let us introduce the set:
$$F_{I,J}=\{(x,y)\in\bR^2: x\in I, y\in
J \; \text{ and the latter inequality holds for} \;  I \; \text{ and} \; J\}.$$

Observe that $\bR^2=\cup_{I,J\in DC(\bR)}F_{I,J}$. Also, $F_{I,J}\sbs I\times
J$. It might be happen that $F_{I_1, J_1}\cap F_{I_2,J_2}\neq 0$
for some different couples of dyadic intervals $(I_1,J_1)$,
$(I_2,J_2)$. Let us take a number $r$ so that $p_+<r<q_-$. Then we
have
\begin{gather*}
\Big\|vM_{\al(\cdot),\bt(\cdot)}^{S,(d)}f\Big\|_{L^{q(\cdot)}(\bR^2)}^r
=\Big\|[vM_{\al(\cdot),\bt(\cdot)}^{S,(d)}f]^r\Big\|_{L^{q(\cdot)/r}(\bR^2)}\leq
  c\sup\limits_{\|h\|_{L^{(q(\cdot)/r)'}(\bR^2)}\leq 1}
  \bigg(\iint_{\bR^2} h[vM_{\al(\cdot),\bt(\cdot)}^{S,(d)}f]^r\bigg).
\end{gather*}

Further, using the above-observed arguments, we find that for such $h$ (we
assume that \linebreak $\|f\wt{M}v\|_{L^{p(\cdot)}(\bR^2)}\leq
1$)
\allowdisplaybreaks
\begin{gather*}
\iint_{\bR^2} h[vM_{\al(\cdot),\bt(\cdot)}^{S,(d)}]^r
\leq \sum_{I,J\in D(\bR)}\iint_{F_{I,J}}h[vM_{\al(\cdot),\bt(\cdot)}^{S,(d)}]^r \\
\leq c\sum_{I,J\in D(\bR)}\bigg(\iint_{I\times J}v^r(x,y)(|I|^{\al(x)}|J|^{\bt(y)})^rh(x,y)dx dy\bigg)\\
\times\bigg(\frac{1}{|I\times J|}\iint_{I\times J} |f(t, \tau)|dt d\tau \bigg)^r \leq c\sum_{I,J\in
D(\bR)}\|(v(\cdot)|I|^{\al(\cdot)}|J|^{\bt(\cdot)})^r
\|_{L^{\frac{q(\cdot)}{r}}(I\times J)}\|h\|_{L^{\big(\frac{q(\cdot)}{r}\big)'}(I\times J)}\\
\times\bigg(\frac{1}{|I| |J|}\iint_{I\times J}|f|\bigg)^r
=c\sum_{I,J\in
D(\bR)}\|(v(\cdot)|I|^{\al(\cdot)}|J|^{\bt(\cdot)})^r\|_{L^{\frac{q(\cdot)}{r}}(I\times J)}\Big(\frac{1}{|I||J|}\iint_{I\times J}|f(t,\tau)|dtd\tau \Big)^r\\
=c\bigg[\sum_{I,J\in
D(\bR)}\|(v(\cdot)|I|^{\al(\cdot)}|J|^{\bt(\cdot)})^r\|_{L^{q(\cdot)}(I\times
J)}
\bigg(\frac{1}{|I| |J|}\iint_{I\times J}|f_1(t, \tau)|dt d\tau \bigg)^r\\
+\sum_{I,J\in
D(\bR)}\|(v(\cdot)|I|^{\al(\cdot)}|J|^{\bt(\cdot)})^r\|_{L^{q(\cdot)}(I\times
J)} \bigg(\frac{1}{|I| |J|}\iint_{I\times J}|f_2(t,\tau)|dt d\tau \bigg)^r\bigg]
=:c[S_1+S_2],
\end{gather*}
where $f_1=f\chi_{\{f\wt{M}_{\al(\cdot),\bt(\cdot)}v\leq 1\}}$,
$f_2=f-f_1$.

Now we estimate $S_1$ and $S_2$ separately. By Corollary A with $\rho \equiv 1$ we have that
\begin{gather*}
S_1=\sum_{I,J\in D(\bR)}(|I||J|)^{-\frac{r}{(p_-)'}} \bigg(\iint_{I\times J}|f_1|(|I| |J|)^{-\frac{1}{p_-}}\|v(\cdot)|I|^{\al(\cdot)}|J|^{\bt(\cdot)}\|_{{}_{L_{I\times J}^{q(\cdot)}}}\bigg)^r\\
=\sum_{I\in D(\bR)}|I|^{-\frac{r}{(p_-)'}}\sum_{J\in
D(\bR)}|J|^{-\frac{r}{(p_-)'}}
\bigg[\int_J\bigg[\int_I|f_1|\Big( \wt{M}_{\al(\cdot),\bt(\cdot)}^{(1)}v\Big)\bigg]\bigg]^r.\\
\end{gather*}
By applying again Corollary A with $\rho\equiv 1$ and generalized
Minkowski inequality, we get
\begin{gather*}
S_{1}\leq c\sum_{I\in D(\bR)}|I|^{-\frac{r}{(p_-)'}}
\bigg(\int_{\bR}\bigg(\int_I|f_1|[\wt{M}_{\al(\cdot), \bt(\cdot)}^{(1)}v]\bigg)^{p_-}\bigg)^{\frac{r}{p_-}}\\
\leq c\sum_{I\in D(\bR)}|I|^{-\frac{r}{(p_-)'}}
\bigg(\bigg(\int_{I}\bigg(\int_{\bR}|f_1|^{p_-}[\wt{M}_{\al(\cdot),
\bt(\cdot)}^{(1)}v]\bigg)^{p_-}\bigg) ^{\frac{1}{p_-}}\bigg)^r\\
\leq c\bigg(\iint_{\bR^2}|f_1|^{p_-}(\wt{M}_{\al(\cdot),\bt(\cdot)}^{(1)}v)\bigg)^{\frac{r}{p_-}}\\
\leq
c\bigg(\iint_{\bR^2}\Big[f(x,y)(\wt{M}_{\al(\cdot),\bt(\cdot)}v)(x,y)\Big]^{p(x,y)}dx dy
\bigg)^{\frac{r}{p_-}}\leq c.
\end{gather*}

By the similar arguments we can see that
$$
S_2\leq
c\bigg(\iint_{\bR^2}[f(x,y)(\wt{M}_{\al(\cdot),\bt(\cdot)}v)(x,y)]^{p(x,y)}dx
dy\bigg)^{\frac{r}{p_+}}\leq c.
$$

Thus we established the desired estimate for the dyadic fractional
maximal function.

Now we pass from $M_{\al(\cdot),\bt(\cdot)}^{S,(d)}$ to
$M_{\al(\cdot),\bt(\cdot)}^S$.

The following inequality for constant $\al$ and $\bt$ was proved
in \cite{KoMeJKMS}  but it is true also for variable $\al$ and
$\bt$:
$$
\Big(M_{a(\cdot), \bt(\cdot)}^{S,(2^k)},f\Big)(x,y)\leq
\frac{C_{\al,\bt}}{|R(0,2^{k+2})|^2}\iint_{R(0,2^{k+2})^2}
S_{t,\tau}(x,y)dt d\tau, \eqno{(4.2)}
$$
where
$$
\Big(M_{\al(\cdot),
\bt(\cdot)}^{S,(2^k)}f\Big)(x,y)=\sup_{\begin{subarray}{c}
I\times J\ni(x,y)\\
|I|,|J|\leq 2^k
\end{subarray}}|I|^{\al(x)-1}|J|^{\bt(y)-1}\iint_{I\times J}|f| , \eqno{(4.3)}$$

$$
S_{t,\tau}(x,y)=\sup_{\begin{subarray}{c}
I-t\ni x\\
J-\tau\ni y\\
I,J\in D(\bR)
\end{subarray}}
|I|^{\al(x)-1}|J|^{\bt(y)-1}\iint_{I-t\times J-\tau} |f|,
\eqno{(4.4)} $$

$$ R(0,r)=\{t:-r\leq t\leq r\}.\notag
$$

Indeed, let $j\in {\Bbb{Z}}$ and let $I$ be an integral such that
$2^{j-1}<|I|\leq 2^j$. Let $j\leq k$, $k\in \bZ$. Suppose that $E$
be the set of those $t\in R(0,2^{k+2})$ for which there is some
$I_1\in D-t$ with $|I_1|=2^{j+1}$ and such that $I\sbs I_1$. Then
(see, e.g., \cite{GaRu}, p. 431)
$$
|E|\geq 2^{k+2},
$$
where $D-t:=\{I-t: I\in D(\bR)\}$.

By the similar arguments, for another interval $J\sbs \bR$, there
is $i\in \bZ$ such that $2^{i-1}<|J|\leq 2^i$. Then for $i\leq k$,
$k\in \bZ$ we have that the set $F$ of those $t\in R(0,2^{k+2})$
for which there is $J_1\in D-t$ with $|J_1|=2^{i+1}$ and $J\sbs
J_1$ has measure greater than or equal to $2^{k+2}$.

To prove (4.2) observe that  for and $(x,y)\in \bR^2$, there are
intervals $Q_1$ and $Q_2$ such that $Q_1\ni x$, $Q_2\ni y$,
$|Q_1|, |Q_2|\leq 2^k$ and
$$
\frac{2}{|Q_1|^{1-\al(x)}|Q_2|^{1-\bt(y)}}\iint_{Q_1\times
Q_2}|f|>(M_{\al(\cdot),\bt(\cdot)}^{S,(2^k)}f)(x,y).
$$

Let $j$ and $i$ be integers such that
$$
2^{j-1}<|Q_1|\leq 2^j,\;\;\;2^{i-1}\leq |Q_2|\leq 2^i.
$$

It is obvious that $j,i\leq k$. Let us define the following sets:
\begin{align*}
E_1:=&\{t\in R(0,2^{k+2}):\exists I\in D(\bR)-t,\\
&|I|=2^{j+1}, Q_1\sbs I\}\\
E_2:=&\{t\in R(0,2^{k+2}):\exists J\in D(\bR)-t,\\
&|J|=2^{j+1}, Q_2\sbs J\}.
\end{align*}

Then using above-observed arguments, we have that $(x\in Q_1\sbs
I, y\in Q_2\sbs J)$:
\begin{gather*}
\frac{1}{2}(M_{\al(\cdot),\bt(\cdot)}^{S,(2^k)}f)(x,y)\leq
\frac{1}{|Q_1|^{1-\al(x)}|Q_2|^{1-\bt(y)}}\iint_{Q_1\times
Q_2}|f|\leq\\
\leq\frac{c_{\al,\bt}}{|I|^{1-\al(x)}|J|^{1-\bt(y)}}\iint_{I\times
J} |f|\leq c_\al S_{t,\tau}(x,y)
\end{gather*}
because $I\in D(\bR)-t$, $J\in  D(\bR)-t$, $I\ni x$, $J\in y$.
Since $|E_1|$, $E_2\geq\frac{|R(0,2^{k+2})|}{2}$, we have that
\begin{gather*}
(M_{\al(\cdot),\bt(\cdot)}^{S,(2^k)}f)(x,y)\leq \frac{c}{|E_1\times E_2|}\iint_{E_1\times E_2} S_{t,\tau}(x,y)dtd\tau\leq\\
\leq\frac{c}{|R(0,2^{k+2})|^2}\iint_{R(0,2^{k+2})^2}S_{t,\tau}(x,y)dtd\tau.
\end{gather*}
Inequality (4.2) is proved.

Further, \allowdisplaybreaks\begin{gather*}
D_{t,\tau}^{(q)}:=\iint_{\bR^2}(S_{t,\tau}(x,y))^{q(x,y)}v(x,y)^{q(x,y)}dxdy=\\
=\iint_{\bR^2}\bigg(\sup_{\begin{subarray}{c}
I-t\ni x\\
J-t\ni y\\
I,J\in D(\bR)
\end{subarray}}
|I|^{\al(x)-1} |J|^{\bt(y)-1}\iint_{(I-t)\times(J-t)}|f(s,\ve)|ds d\ve\bigg)^{q(x,y)} v(x,y)^{q(x,y)}dx dy=\\
=\iint_{\bR^2} \bigg(\sup_{\begin{subarray}{c}
I-t\ni x\\
J-t\ni y\\
I,J\in D(\bR)
\end{subarray}}
|I|^{\al(x-t)-1} |J|^{\bt(y-\tau)-1}\iint_{I\times J}f(s-t,s-t)ds d\ve\bigg)^{q(x-t,y-\tau)}\times\\
\times v(x-t,y-\tau)^{q(x-t,y-\tau)}dxdy=\\
=\iint_{\bR^2}(M_{\al(\cdot-t),\bt(\cdot-\tau)}^{S,(d)} f(\cdot-t,
\cdot-\tau))^{q(x-t,y-\tau)}v^{q(x-t,y-\tau)}(x-t,y-\tau)dx dy.
\end{gather*}

Observe now that
$$\frac{1}{p(\cdot-t,\cdot-\tau)}-\frac{1}{(q(\cdot-t,\cdot-\tau))_+}<(\al(\cdot-t,\cdot-\tau))_-$$

$$
\leq (\al(\cdot-t,\cdot
-\tau))_+<\frac{1}{(p(\cdot-t,\cdot-\tau))_-}$$ and since the
constants in the estimates of
$\|(M_{\al(\cdot),\bt(\cdot)}^{S,(d)}
f)\|_{{}_{L_{v(\bR^2)}^{q(\cdot)}}}$
depends only on  $p_+,\; q_-,$ we have that

$$ D_{t,\tau}^{(q)}\leq c $$
because
$$
I:=\iint_{\bR^2}|f(x-t,t-\tau)|^{p(x-t,y-\tau)}(\wt{M}_{\al(\cdot-t),\beta(\cdot
-\tau)}v(\cdot -t,\cdot -\tau))^{p(x-t,y-\tau)} (x,y)dx dy$$

$$
=\iint_{\bR^2}|f(x,y)|^{p(x,y)}(\wt{M}_{\al(\cdot),\bt(\cdot)}v)^{p(x,y)}(x,y)
dxdy\leq 1.\notag \eqno{(4.5)}
$$
Now let us see that equality (5) holds.

First observe that
\begin{gather*}
(\wt{M}_{\al(\cdot -t),\bt(\cdot -\tau)}^{(1)}v(\cdot -t,\cdot
-\tau))(x,y)=\sup_{\begin{subarray}{c}
I\ni x\\
J\ni y
\end{subarray}}|I\times J|^{-\frac{1}{(p(\cdot -t,\cdot -\tau))_-}}\\
\times \|v(\cdot -t, \cdot -\tau)|I|^{\al(\cdot -t)}|\tau|^{\bt(\cdot -\tau)}
\|_{L^{q(\cdot -t,\cdot -\tau)}(I\times J)}=\\
=\sup_{\begin{subarray}{c}
I\ni x\\
J\ni y
\end{subarray}}|I\times J|^{-\frac{1}{p_-}} \inf\bigg\{\lb>0:\!\!\iint_{I\times J}\!\!\Big[\frac{v(x-t, y-\tau)|I|^{\al(x-t)}|J|^{\bt(y-\tau)}}{\lb}\Big]^{q(x-t,y-\tau)}dx dy\!\leq \!1 \bigg\}\!=\\
=\sup_{\begin{subarray}{c}
I-t\ni x-t\\
J-\tau\ni y-\tau
\end{subarray}}|I-t|^{-\frac{1}{p_-}}|J-\tau|^{\frac{1}{p_-}}\inf\bigg\{\lb>0:\!\!\!\iint_{(I-t)\times(J-\tau)}\!\!\!\Big[\frac{v(x,y)|I|^{\al(x)}
|J|^{\bt(y)}}{\lb}\Big]^{q(x,y)}dx dy\!\leq\! 1\bigg\}\!=\\
=(\wt{M}_{\al(\cdot),\bt(\cdot)}^{(1)}v)(x-t,y-\tau).
\end{gather*}

Analogous estimates hold also for $\bigg(\wt{M}^{(2)}_{\al(\cdot -t),\bt(\cdot -\tau)}v\bigg) (\cdot
-t, \cdot -\tau)$.

Hence
\begin{gather*}
I=\iint_{\bR^2}|f(x-t,t-\tau)|^{p(x-t,y-\tau)}(\wt{M}_{\al(\cdot),\bt(\cdot)}v)^{p(x-t,y-\tau)}
(x-t,y-\tau)dx dy=\\
=\iint_{\bR^2}|f(x,y)|^{p(x,y)}(\wt{M}_{\al(\cdot),\bt(\cdot)}v)^{p(x,y)}(x,y)
dxdy.
\end{gather*}

Thus, we have seen that
\begin{gather*}
D_{t,\tau}^{(q)}\leq c\;\;\;\text{if}\;\;\;\|f(\wt{M}_{\al(\cdot),\bt(\cdot)}v)\|_{L^{p(\cdot)}(\bR^2)}\leq 1\Longleftrightarrow\\
\Longleftrightarrow \|S_{t,\tau}(\cdot,\cdot)v(\cdot,\cdot)\|_{L^{q(\cdot)}(\bR^2)}\leq c\|f(\wt{M}_{\al(\cdot),\bt(\cdot)}v)\|_{L^{p(\cdot)}(\bR^2)}\Longleftrightarrow \\
\Longleftrightarrow \sup_{\|h\|_{L^{q(\cdot)}(\bR^2)}\leq 1}\iint_{\bR^2} S_{t,\tau}(x,y)v(x,y)h(x,y)dx dy\leq c
\|f(\wt{M}_{\al(\cdot),\bt(\cdot)}v)\|_{L^{p(\cdot,\cdot)}(\bR^2)}.
\end{gather*}

For $h$, $\|h\|_{L^{q'(\cdot)}(\bR^2)}\leq 1$, we have
that (recall that (4.2) holds):
\begin{gather*}
\iint_{\bR^2}(M_{\al(\cdot), \bt(\cdot)}^{S,(2^k)},f)(x,y) v(x,y)h(x,y)dx dy\leq \\
\leq c|R(0,2^{k+2})|^{-2}\iint_{\bR^2}\bigg[R(0,2^{k+2})S_{t,\tau}(x,y)dt d\tau\bigg]v(x,y) h(x,y)dx dy=\\
= c|R(0,2^{k+2})|^{-2}\iint_{R(0,2^{k+2})^2}\bigg(\iint_{\bR^2}S_{t,\tau(x,y)}v(x,y)h(x,y)dx dy\bigg)dt d\tau\leq\\
\leq c\|f
(\wt{M}_{\al(\cdot,\bt(\cdot))}v)\|_{L^{p(\cdot,\cdot)}(\bR^2)}|R(0,2^{k+2})|^{-2}
\iint_{R(0,2^{k+2})^2}dtd\tau=c\|f(\wt{M}_{\al(\cdot),\bt(\cdot)}v)\|_{L^{p(\cdot,\cdot)}(\bR^2)}.
\end{gather*}
Passing now $k$ to the infinity and taking the supremum with
respect to $h$ in the last inequality, we get the desired result.
\vskip+0.2cm

{\em Remark} 4.3.  Observe that for $p\equiv \const$, $\al\equiv
\const$, $\bt\equiv \const$ and $v\in L^{q(\cdot,\cdot)}(\bR^2)$,
the estimate
$$
(\wt{M}_{\al(\cdot), \bt(\cdot)}v)(x,y)\leq
\sup_{\begin{subarray}{c}
I\ni x\\
J\ni y
\end{subarray}}|I|^{\al-\frac{1}{p}}|J|^{\bt-\frac{1}{p}}\bigg(\iint_{I\times J}v(x,y)^{q(x,y)}dx dy\bigg)^{\frac{1}{q_+}}=: (\ol{M}_{\alpha, \beta}v)(x,y)
$$
holds. \vskip+0.2cm

{\bf Corollary 4.5.} {\em  Let $p$, $\al$ and $\bt$ be constant
and let $\frac{1}{p}-\frac{1}{q_+}<\al$, $\bt<\frac{1}{p}$.
Suppose that $v\in L^{q(\cdot,\cdot)}(\bR^2)$. Then the
inequality
$$
\|v(M_{\al(\cdot),\bt(\cdot)}f)\|_{L^{q(\cdot,\cdot)}(\bR^2)}\leq
c\|f(\cdot,\cdot)(\ol{M}_{'\al,\bt})(\cdot,\cdot) \|_{L^p(\bR^2)}
$$
holds.}

\vskip+0.2cm

 Corollary 4.5 follows immediately from Theorem
4.1  and Remark 4.3. \vskip+0.2cm

 {\em Proof of Corollary}  $4.2$. This proposition will be proved if we
show that $(\wt{M}_{\al(\cdot),\bt(\cdot)}v)(x,y)\leq c$ in
Theorem 4.1. Indeed, if the condition

$$
A:=\sup_{I,J\sbs \bR}
\||I|^{\al(\cdot)}|J|^{\bt(\cdot)}v(\cdot)\|_{L^{q(\cdot)}(I\times
J)}|I\times J|^{-\frac{1}{{\ol{p}}_{I\times J}}}<\infty,
$$
is satisfied, where
$$
\ol{p}_{I\times J}=\begin{cases}
p_-,\;\;\;\text{if}\;\;\;|I||J|\leq 1,\\
p_+,\;\;\;\text{if}\;\;\;|I||J|> 1,
\end{cases}
$$
then
$$
\||I|^{\al(\cdot)}|J|^{\bt(\cdot)}v(\cdot)\|_{L^{q(\cdot)}(I\times
J)} |I\times J|^{-\frac{1}{p_+}}\leq A<\infty
$$
and
$$
\hskip+3.5cm\||I|^{\al(\cdot)}|J|^{\bt(\cdot)}v(\cdot)\|_{L^{q(\cdot)}(I\times
J)} |I\times J|^{-\frac{1}{p_-}}\leq A<\infty.
\hskip+3.5cm\square
$$


\noindent\emph{Proof of Theorem $4.3$}. Let us recall that by the
symbol $M_{\al(\cdot),\bt(\cdot)}^{S,(d)}$ is denoted the dyadic
strong fractional maximal operator.

\emph{Sufficiency.} We use the notation of the proof of Theorem
4.1. First we construct the sets $F_{I\times J}$.

Take $r$ so that $p<r<q_-$ and observe that
\begin{gather*}
\|v(M_{\al(\cdot),\bt(\cdot)}^{S,(d)}f)\|_{L^{q(\cdot)}({\Bbb{R}}^2)}\leq\\
\leq c\sup\limits_{\|h\|_{L^{(q(\cdot)/r)'}({\Bbb{R}}^2)}\leq
1}\bigg(\iint_{\bR^2}h[v M_{\al(\cdot),\bt(\cdot)}^{s,(d)}
f]^r\bigg).
\end{gather*}

Let $\|f\|_{L^{p}_{w}(\mathbb{R}^{2})}\leq 1$. Then for such an
$h$ we have that

\allowdisplaybreaks\begin{gather*} S=\iint_{\bR^2}h[v
M_{\al(\cdot),\bt(\cdot)}^{S,(d)}f]^r
\leq\sum_{I,J\in D(\bR)}\iint_{F_{I,J}}h[v M_{\al(\cdot),\bt(\cdot)}^{s,(d)}]^r\\
\leq c \sum_{I,J\in D(\bR)}\bigg(\iint_{I\times J}
v^r(|I||^{\al(x)}|J|^{\bt(y)})^r h(x,y)dx dy\bigg)
\Big(\frac{1}{|I||J|}\iint_{I\times J}|f(t, \tau)|dt d\tau \Big)^r\\
\leq c\sum_{I,J\in D(\bR)}
\|(v(\cdot)|I|^{\al(\cdot)}|J|^{\bt(\cdot)})^r\|_{L^{q(\cdot)/r}
(I\times J)}|h|_{L^{(q(\cdot)/r)'}(I\times J)}
\Big(\frac{1}{|I||J|}\iint_{I\times J}|f(t,\tau)|dt d\tau \Big)^r\\
=c\sum_{I,J\in
D(\bR)}\|v(\cdot)|I|^{\al(\cdot)}|J|^{\bt(\cdot)}\|_{L^{q(\cdot)}
(I\times J)}^r \Big(\frac{1}{|I||J|}\iint_{I\times J}|f(t, \tau)|dt d\tau\Big)^{r}.
\end{gather*}
By the condition of theorem we get that
\begin{gather*}
S\leq c\sum_{I,J\in D(\bR)}\bigg(\int_I
w_1^{-p'}\bigg)^{-\frac{r}{p'}}\bigg(\int_J
w_2^{-p'}\bigg)^{-\frac{r}{p'}}\bigg(\iint_{I\times J}|f|\bigg)^r
\end{gather*}
Applying Corollary A with $\rho \equiv 1$ we
derive the following estimates:
\begin{gather*}
S\leq c\sum_{J\in D(\bR)}\bigg(\int_J w_2^{-p'}\bigg)^{-\frac{r}{p'}}\bigg(\iint_{\bR}w_1(t)^p\bigg(\int_J|f(t,\tau)|d\tau\bigg)^pdt\bigg)^{\frac{r}{p}}\\
\leq c\sum_{J\in D(\bR)}\bigg(\int_J w_2^{-p'}\bigg)^{-\frac{r}{p'}}\bigg(\int_J\bigg(\int_{\bR}w_1^p(t)|f(t,\tau)|^pdt\bigg)^{\frac{1}{p}}d\tau\bigg)^r\\
\leq c\bigg(\iint_{\bR^2}|f(t,\tau)|^pw^p(t,\tau)dt
d\tau\bigg)^{r/p}\leq c.
\end{gather*}

Thus we established the desired inequality for the dyadic fractional
maximal function.

Now we can pass to the fractional maximal function
$M_{\al(\cdot),\bt(\cdot)}^S$ in the same manner as in the proof
of Theorem 4.1.

\emph{Necessity} follows easily by taking appropriate test functions in the two--weight inequality. Details  are omitted. \hfill $\square$ \vskip+0.2cm

\noindent\emph{Proof of Corollary $4.3$.} The proof is a direct
consequence of Theorem 4.3 and the fact that the condition $p\in
\mathcal{P}(\bR^2)$ implies the inequality

$$
\|f w\|_{L^{p_-}(I\times J)}\leq c\|fw\|_{L^{p(\cdot)}(\bR^2)}.
\;\;\;\;\;\; \Box
$$

Corollary 4.4 follows from the fact: $p\in
\mathcal{P}_\infty(\bR^2)\Longrightarrow p\in \mathcal{P}(\bR^2)$,
provided that $p_-=p(\infty)$, and Corollary 4.3.

\

\subsection*{The Case of a Bounded Domain}

\

Let $\Om$ be a bounded domain in $\bR^2$ and let
$$
(M_{\al(\cdot),\bt(\cdot)}^{S,R_0}f)(x,y)=\sup_{\begin{subarray}{c}
R=I\times J\\
R\ni(x,y)
\end{subarray}}|I|^{\al(x)-1}|J|^{\bt(y)-1}\iint_{R}|f(t, \tau)| dt d\tau
$$
where $(x,y)\in\Om$. For simplicity assume that $\Om=R_0$, where
$R_0$ is a fixed rectangle in $\bR^2$.

Taking the results of the previous subsections into account  we
can formulate the following statements proofs of which are
omitted:

\vskip+0.2cm

{\bf Theorem 4.4 (Fefferman-Stein type inequality).} {\em Let
$1<p_-(R_0)\leq p_+(R_0)<q_-(R_0)\leq q_+(R_0)<\infty$ and let
\begin{gather*}
\frac{1}{p_-(R_0)}-\frac{1}{q_-(R_0)}<\al_-(R_0)\leq \al_+(R_0)<\frac{1}{p_-(R_0)},\\
\frac{1}{p_-(R_0)}-\frac{1}{q_-(R_)}<\bt_-(R_0)\leq\al_+(R_0)<\frac{1}{p_-(R_0)}.
\end{gather*}
Then there is a positive constant $b$ such that the following
inequality
$$
\|(M_{\al(\cdot),\bt(\cdot)}^{S,R_0}f)v\|_{L^{q(\cdot)}(R_0)}\leq
b\|f(\wt{M}_{\al(\cdot),\bt(\cdot)}^{R_0}v)\|_{L^{p(\cdot)}(R_0)},
$$
holds, where
$$
(\wt{M}_{\al(\cdot),\bt(\cdot)}^{(R_0)}v)(x,y)=\sup_{\begin{subarray}{c}
R=I\times J\\
R\ni (x,y)
\end{subarray}}
|R|^{-\frac{1}{p_-}}\|v(\cdot)|I|^{\al(\cdot)}|J|^{\bt(\cdot)}\|_{L^{q(\cdot)}(R\cap
R_0)},\;\;\;(x,y)\in R_0.
$$}

\vskip+0.1cm

{\em Remark} 4.4.  Let $\al(x)\equiv\al\equiv \const$,
$\bt(x)\equiv\bt\equiv\const$ in Theorem 4.4. Then it is easy to
see that
$$
(\wt{M}_{\al,\bt}^{(R_0)}v)(x,y)\leq \sup_{\begin{subarray}{c}
R=I\times J\\
R\ni (x,y)
\end{subarray}}
|I|^{\al-\frac{1}{p_-}}|J|^{\bt-\frac{1}{p_-}}\bigg(\iint_{R\cap
R_0}v^{q(x,y)}(x,y)dxdy\bigg)^{\frac{1}{q_+}}.
$$
\vskip+0.2cm

{\bf Theorem 4.5 (Trace inequality).} {\em Let $1<p_-(R_0)\leq
p_+(R_0)<q_-(R_0)\leq q_+(R_0)<\infty$ and let
$\frac{1}{p_-(R_0)}-\frac{1}{q_+(R_0)}<\al_-(R_0)\leq
\al_+(R_0)<\frac{1}{p_-(R_0)}$,
$\frac{1}{p_-(R_0)}-\frac{1}{q_+(R_0)}<\bt_-(R_0)\leq\bt_+(R_0)<\frac{1}{p_-(R_0)}$.
Suppose that the weight function $v$ on $\bR^2$ satisfies the
condition
$$
\sup_{\begin{subarray}{c}
R=I\times J\\
R\sbs R_0
\end{subarray}} \||I|^{\al(\cdot)}|J|^{\bt(\cdot)}v(\cdot)\|_{L^{q(\cdot)}({\Bbb{R}})}|R|^{-\frac{1}{p_-(R)}}<\infty.
$$
Then $M_{\al(\cdot),\bt(\cdot}^{S,R_0}$ is bounded from
$L^{p(\cdot)}(R_0)$ to $L_v^{q(\cdot)}(R_0)$.} \vskip+0.2cm

{\bf Theorem 4.6.} {\em  Let $R_0:=I_0\times J_0$,
$1<p_-(R_0)<q_-(R_0)\leq q_+(R_0)<\infty$ with $p_+(R_0)<\infty$,
$0<\al_-(R_0)\leq \al_+(R_0)<1$ and $0<\bt_-(R_0)\leq
\bt_+(R_0)<1$. Suppose that $v$ and $w$ are weights on $R_0$ with
$w(x,y)=w_1(x)w_2(y)$, where $w_1\in RD^{(d)}(I_0)$, $w_2\in
RD^{(d)}(J_0)$. If
$$
\sup_{\begin{subarray}{c}
R\sbs R_0\\
R=I\times J
\end{subarray}} |R|^{-1}\|v(\cdot)|I|^{\al(\cdot)}|J|^{\bt(\cdot)}\|_{L^{q(\cdot)}(R)}\|w^{-1}\|_{L^{(p_-)}(R)}<\infty,
$$
then $M_{\al(\cdot),\bt(\cdot)}^{S,R_0}$ is bounded from
$L_w^{p(\cdot)}(R_0)$ to $L_v^{p(\cdot)}(R_0)$.}

\section*{Acknowledgements}
The  authors were partially supported by the Georgian National
Science Foundation Grant (project numbers: No. GNSF/ST09/23/3-100
and No. GNSF/ST07/3-169). \vskip+1cm

\vskip+0.5cm

 Authors' Addresses:

\

V. Kokilashvili:

\ A. Razmadze Mathematical Institute, 1. M. Aleksidze Str., 0193
Tbilisi, Georgia

Second Address:  Faculty of Exact and Natural Sciences, I.
Javakhishvili Tbilisi State University 2, University St., Tbilisi
0143 Georgia.

E-mail: kokil@@rmi.acnet.ge

\vskip+0.5cm

A. Meskhi: \

A. Razmadze Mathematical Institute, M. Aleksidze St.,  Tbilisi
0193, Georgia

Second Address: Department of Mathematics,  Faculty of Informatics
and Control Systems, Georgian Technical University, 77, Kostava
St., Tbilisi, Georgia.

e-mail:  meskhi@@rmi.acnet.ge \vskip+0.5cm

\end{document}